\documentclass[12pt,a4paper,reqno]{amsart}
\usepackage{amsmath}
\usepackage{amsfonts}
\usepackage{amssymb}
\numberwithin{equation}{section}

\theoremstyle{plain}

\newtheorem{theorem}[subsection]{Theorem}
\newtheorem{proposition}[subsection]{Proposition}
\newtheorem{lemma}[subsection]{Lemma}

\newtheorem*{inductive}{Lemma \ref{ind-step}}

\theoremstyle{definition}

\newtheorem{definition}[subsection]{Definition}

\theoremstyle{remark}

\newtheorem*{example}{Example}

\newtheorem*{remark}{Remark}
\newtheorem*{remarks}{Remarks}

\renewcommand{\leq}{\leqslant}
\renewcommand{\geq}{\geqslant}

\newsavebox{\proofbox}
\savebox{\proofbox}{\begin{picture}(7,7)%
  \put(0,0){\framebox(7,7){}}\end{picture}}

\def\proof{\noindent\textit{Proof. }}
\def\endproof{\hfill{\usebox{\proofbox}}}
\def\E{\mathbb{E}}

\def\F{\mathbb{F}}
\def\Z{\mathbb{Z}}
\def\R{\mathbb{R}}
\def\T{\mathbb{T}}
\def\C{\mathbb{C}}

\def\M{\mathbf{M}}
\def\S{\mathcal{S}}
\def\eps{\varepsilon}
\newcommand\Spec{\operatorname{Spec}}
\newcommand\Bohr{\operatorname{Bohr}}

\newcommand\Supp{\operatorname{Supp}}

\begin{document}

\title[Quantitative idempotent theorem]{A quantitative version of the idempotent theorem in harmonic analysis}

\author{Ben Green}
\address{Centre for Mathematical Sciences, Wilberforce Road, Cambridge CB3 0WA, England.}
\email{b.j.green@dpmms.cam.ac.uk}

\author{Tom Sanders}
\address{Centre for Mathematical Sciences, Wilberforce Road, Cambridge CB3 0WA, England.
}
\email{t.sanders@dpmms.cam.ac.uk}

\thanks{The first author is a Clay Research Fellow, and is pleased to acknowledge the support of the Clay Mathematics Institute.}

\begin{abstract} 
Suppose that $G$ is a locally compact abelian group, and write $\M(G)$ for the algebra of bounded, regular, complex-valued measures under convolution. A measure $\mu \in \M(G)$ is said to be \emph{idempotent} if $\mu \ast \mu = \mu$, or alternatively if $\widehat{\mu}$ takes only the values $0$ and $1$.  The Cohen-Helson-Rudin idempotent theorem states that a measure $\mu$ is idempotent if and only if the set $\{\gamma \in \widehat{G} : \widehat{\mu}(\gamma) = 1\}$ belongs to the \emph{coset ring} of $\widehat{G}$, that is to say we may write
\[ \widehat{\mu}  = \sum_{j = 1}^L \pm 1_{\gamma_j + \Gamma_j}\]
where the $\Gamma_j$ are open subgroups of $\widehat{G}$.

In this paper we show that $L$ can be bounded in terms of the norm $\Vert \mu \Vert$, and in fact one may take $L \leq \exp\exp(C \Vert \mu \Vert^4)$. In particular our result is non-trivial even for finite groups.
\end{abstract}

\maketitle

\section{Introduction}

Let us begin by stating the idempotent theorem. Let $G$ be a locally compact abelian group with dual group $\widehat{G}$. Let $\M(G)$ denote the \emph{measure algebra} of $G$, that is to say the algebra of bounded, regular, complex-valued measures on $G$. We will not dwell on the precise definitions here since our paper will be chiefly concerned with the case $G$ finite, in which case $\M(G) = L^1(G)$. For those parts of our paper concerning groups which are not finite, the book \cite{rudin-book} may be consulted. A discussion of the basic properties of $\M(G)$ may be found in Appendix E of that book.

If $\mu \in \M(G)$ satisfies $\mu \ast \mu = \mu$, we say that $\mu$ is \emph{idempotent}. Equivalently, the Fourier-Stieltjes transform $\widehat{\mu}$ satisfies $\widehat{\mu}^2 = \widehat{\mu}$ and is thus $0,1$-valued. 
\begin{theorem}[Cohen's idempotent theorem] $\mu$ is idempotent if and only if $\{\gamma \in \widehat{G} : \widehat{\mu}(\gamma) = 1\}$ lies in the coset ring of $\widehat{G}$, that is to say 
\[ \widehat{\mu} = \sum_{j = 1}^L \pm 1_{\gamma_j + \Gamma_j},\]
where the $\Gamma_j$ are open subgroups of $\widehat{G}$.
\end{theorem}

This result was proved by Paul Cohen \cite{cohen}. Earlier results had been obtained in the case $G = \T$ by Helson \cite{helson} and $G = \T^d$ by Rudin \cite{rudin-paper}. See \cite[Ch. 3]{rudin-book} for a complete discussion of the theorem.

When $G$ is finite the idempotent theorem gives us no information, since $\M(G)$ consists of \emph{all} functions on $G$, as does the coset ring. The purpose of this paper is to prove a quantitative version of the idempotent theorem which does have non-trivial content for finite groups. 

\begin{theorem}[Quantitative idempotent theorem]\label{mainthm} Suppose that $\mu \in \M(G)$ is idempotent. Then we may write
\[ \widehat{\mu} = \sum_{j = 1}^L \pm 1_{\gamma_j + \Gamma_j},\] where $\gamma_j \in \widehat{G}$, each $\Gamma_j$ is an open subgroup of $\widehat{G}$ and $L \leq e^{e^{C\Vert \mu \Vert^4}}$ for some absolute constant $C$. The number of distinct subgroups $\Gamma_j$ may be bounded above by $\Vert \mu \Vert + \frac{1}{100}$.
\end{theorem}
\begin{remark} In this theorem (and in Theorem \ref{mainthm-finite} below) the bound of $\Vert \mu \Vert + \frac{1}{100}$ on the number of different subgroups $\Gamma_j$ (resp. $H_j$) could be improved to $\Vert \mu \Vert + \delta$, for any fixed positive $\delta$. We have not bothered to state this improvement as obtaining the correct dependence on $\delta$ would add unnecessary complication to an already technical argument. Furthermore the improvement is only of any relevance at all when $\Vert \mu \Vert$ is a tiny bit less than an integer.\end{remark}

To apply Theorem \ref{mainthm} to finite groups it is natural to switch the r\^oles of $G$ and $\widehat{G}$. One might also write $\widehat{\mu} = f$, in which case the idempotence of $\mu$ is equivalent to asking that $f$ is $0,1$-valued, or the characteristic function of a set $A \subseteq G$. It turns out to be just as easy to deal with functions which are $\Z$-valued. The norm $\Vert \mu \Vert$ is the $\ell^1$-norm of the Fourier transform of $f$, also known as the \emph{algebra norm} $\Vert f \Vert_A$ or sometimes, in the computer science literature, as the \emph{spectral norm}. We will define all of these terms properly in the next section.

\begin{theorem}[Main theorem, finite version]\label{mainthm-finite}
Suppose that $G$ is a finite abelian group and that $f : G \rightarrow \Z$ is a function with $\Vert f \Vert_A \leq M$. Then we may write
\[ f = \sum_{j = 1}^L \pm 1_{x_j + H_j},\] where $x_j \in G$, each $H_j \leq G$ is a subgroup and $L \leq e^{e^{CM^4}}$. Furthermore the number of distinct subgroups $H_j$ may be bounded above by $M + \frac{1}{100}$.
\end{theorem}

Theorem \ref{mainthm-finite} is really the main result of this paper. Theorem \ref{mainthm} is actually deduced from it (and the ``qualititative'' version of the idempotent theorem). This reduction is contained in Appendix \ref{lcag-finite}. The rest of the paper is entirely finite in nature and may be read independently of Appendix \ref{lcag-finite}.

\section{Notation and conventions}

Background for much of the material in this paper may be found in the book of Tao and Vu \cite{tao-vu-book}. We shall often give appropriate references to that book as well as the original references. Part of the reason for this is that we hope the \emph{notation} of \cite{tao-vu-book} will become standard. 

\textsc{constants.} Throughout the paper the letters $c,C$ will denote absolute constants which could be specified explicitly if desired. These constants will generally satisfy $0 < c \ll 1 \ll C$. Different instances of the notation, even on the same line, will typically denote different constants. Occasionally we will want to fix a constant for the duration of an argument; such constants will be subscripted as $C_0,C_1$ and so on. 

\textsc{measures on groups.} Apart from in Appendix \ref{lcag-finite} we will be working with functions defined on finite abelian groups $G$. As usual we write $\widehat{G}$ for the group of characters $\gamma : G \rightarrow \C^{\times}$ on $G$. We shall always use the \emph{normalised counting measure} on $G$ which attaches weight $1/|G|$ to each point $x \in G$, and \emph{counting measure} on $\widehat{G}$ which attaches weight one to each character $\gamma \in \widehat{G}$. Integration with respect to these measures will be denoted by $\E_{x \in G}$ and $\sum_{\gamma \in \widehat{G}}$ respectively. Thus if $f : G \rightarrow \C$ is a function we define the $L^p$-norm
\[ \Vert f \Vert_p := \big( \E_{x \in G} |f(x)|^p \big)^{1/p} = \big( \frac{1}{|G|} \sum_{x \in G} |f(x)|^p \big)^{1/p},\] whilst the $\ell^p$-norm of a function $g : \widehat{G} \rightarrow \C$ is defined by
\[ \Vert g \Vert_p := \big( \sum_{\gamma \in \widehat{G}} |g(\gamma)|^p \big)^{1/p}.\]
The group that any given function is defined on will always be clear from context, and so this notation should be unambiguous.

\textsc{fourier analysis.} If $f : G \rightarrow \C$ is a function and $\gamma \in \widehat{G}$ we define the Fourier transform $\widehat{f}(\gamma)$ by
\[ \widehat{f}(\gamma) := \E_{x \in G} f(x) \overline{\gamma(x)}.\] We shall sometimes write this as $(f)^{\wedge}(\gamma)$ when $f$ is given by a complicated expression.
If $f_1,f_2 : G \rightarrow \C$ are two functions we define their convolution by
\[ f_1 \ast f_2(t) := \E_{x \in G} f_1(x) f_2(t - x).\]
We note the basic formul{\ae} of Fourier analysis:
\begin{enumerate}
\item (Plancherel) $\langle f_1 , f_2\rangle := \E_{x \in G} f_1(x) \overline{f_2(x)} = \sum_{\gamma \in \widehat{G}} \widehat{f}_1(\gamma) \overline{\widehat{f}_2(\gamma)} = \langle \widehat{f}_1, \widehat{f}_2\rangle$;
\item (Inversion) $f(x) = \sum_{\gamma \in \widehat{G}} \widehat{f}(\gamma) \gamma(x)$;
\item (Convolution) $(f_1 \ast f_2)^{\wedge} = \widehat{f}_1 \widehat{f}_2$.
\end{enumerate}
In this paper we shall be particularly concerned with the algebra norm 
\[ \Vert f \Vert_A := \Vert \widehat{f}\Vert_1 = \sum_{\gamma \in \widehat{G}} |\widehat{f}(\gamma)|.\]
The name comes from the fact that it satisfies $\Vert f_1f_2 \Vert_1 \leq \Vert f_1 \Vert_A \Vert f_2 \Vert_A$ for any $f_1,f_2 : G \rightarrow \C$.

If $f : G \rightarrow \C$ is a function then we have $\Vert \widehat{f} \Vert_{\infty} \leq \Vert f \Vert_1$ (a simple instance of the Hausdorff-Young inequality). If $\rho \in [0,1]$ is a parameter we define
\[ \Spec_{\rho}(f) := \{ \gamma \in \widehat{G} : |\widehat{f}(\gamma)| \geq \rho\Vert f \Vert_1\}.\]

\textsc{freiman isomorphism.} Suppose that $A \subseteq G$ and $A' \subseteq G'$ are subsets of abelian groups, and that $s \geq 2$ is an integer. We say that a map $\phi : A \rightarrow A'$ is a Freiman $s$-homomorphism if $a_1 + \dots + a_s = a_{s+1} + \dots + a_{2s}$ implies that $\phi(a_1) + \dots + \phi(a_s) = \phi(a_{s+1}) + \dots + \phi(a_{2s})$. If $\phi$ has an inverse which is also a Freiman $s$-homomorphism then we say that $\phi$ is a Freiman $s$-isomorphism and write $A \cong_s A'$. 

\section{The main argument}\label{main-deduction-sec}

In this section we derive Theorem \ref{mainthm-finite} from Lemma \ref{ind-step} below. The proof of this lemma forms the heart of the paper and will occupy the next five sections. 

Our argument essentially proceeds by induction on $\Vert f \Vert_A$, splitting $f$ into a sum $f_1 + f_2$ of two functions and then handling those using the inductive hypothesis. As in our earlier paper \cite{green-sanders}, it is not possible to effect such a procedure entirely within the ``category'' of $\Z$-valued functions. One must consider, more generally, functions which are $\eps$-almost $\Z$-valued, that is to say take values in $\Z + [-\eps,\eps]$. If a function has this property we will write $d(f,\Z) < \eps$. In our argument we will always have $\eps < 1/2$, in which case we may unambiguously define $f_{\Z}$ to be the integer-valued function which most closely approximates $f$.

\begin{lemma}[Inductive Step]\label{ind-step}
Suppose that $f : G \rightarrow \R$ has $\Vert f \Vert_A \leq M$, where we take $M \geq 1$, and that $d(f,\Z) \leq e^{-C_1M^4}$. Set $\eps := e^{-C_0M^4}$, for some constant $C_0$.
Then we may write $f = f_1 + f_2$, where 
\begin{enumerate}
\item either $\Vert f_1 \Vert_A \leq \Vert f \Vert_A - 1/2$ or else $(f_1)_{\Z}$ may be written as $\sum_{j = 1}^L \pm 1_{x_j + H}$, where $H$ is a subgroup of $G$ and $L \leq e^{e^{C'(C_0)M^4}}$;
\item $\Vert f_2 \Vert_A \leq \Vert f \Vert_A - \frac{1}{2}$ and 
\item $d(f_1,\Z) \leq d(f,\Z) + \eps$ and $d(f_2,\Z) \leq 2d(f,\Z) + \eps$.
\end{enumerate}
\end{lemma}

\emph{Proof of Theorem \ref{mainthm-finite} assuming Lemma \ref{ind-step}.} We apply Lemma \ref{ind-step} iteratively, starting with the observation that if $f : G \rightarrow \Z$ is a function then $d(f,\Z) = 0$. Let $\eps = e^{-C_0M^4}$ be a small parameter, where $C_0$ is much larger than the constant $C_1$ appearing in the statement of Lemma \ref{ind-step}. Split 
\[ f = f_1 + f_2\]
according to Lemma \ref{ind-step} in such a way that $d(f_1,\Z),d(f_2,\Z) \leq \eps$. Each $(f_i)_{\Z}$ is a sum of at most $e^{e^{CM^4}}$ functions of the form $\pm 1_{x_j + H_i}$ (in which case we say it is \emph{finished}), or else we have $\Vert f_i \Vert_A \leq \Vert f \Vert_A - \frac{1}{2}$. 

Now split any unfinished functions using Lemma \ref{ind-step} again, and so on (we will discuss the admissibility of this shortly). After at most $2M - 1$ steps all functions will be finished. Thus we will have a decomposition
\[ f = \sum_{k = 1}^L f_k,\]
where 
\begin{enumerate}
\item[(a)] $L \leq 2^{2M-1}$;
\item[(b)] for each $k$, $(f_k)_{\Z}$ may be written as the sum of at most $e^{e^{CM^4}}$ functions of the form $\pm 1_{x_{j,k} + H_k}$, where $H_k \leq G$ is a subgroup, and
\item[(c)] $d(f_k,\Z) \leq 2^{2M}\eps$ for all $k$.
\end{enumerate} 
The last fact follows by an easy induction, noting carefully the factor of 2 in (iii) of Lemma \ref{ind-step}. Note that as a consequence of this, and the fact that $C_0 \ggg C_1$, our repeated applications of Lemma \ref{ind-step} were indeed valid.

Now we clearly have
\[ \Vert f - \sum_{k = 1}^L (f_k)_{\Z} \Vert_{\infty} \leq 2^{4M-1}\eps < 1.\] Since $f$ is $\Z$-valued we are forced to conclude that in fact
\[ f = \sum_{k = 1}^L (f_k)_{\Z}.\]
It remains to establish the claim that we may take $L \leq M + \frac{1}{100}$.
By construction we have 
\[ \Vert f \Vert_A = \sum_{k = 1}^L \Vert f_k \Vert_A.\]
If $(f_k)_{\Z}$ is not identically 0 then, since $\eps$ is so small, we have from (c) above that
\[ \Vert f_k \Vert_A \geq \Vert f_k \Vert_{\infty} \geq \Vert (f_k)_{\Z} \Vert_{\infty} - 2^{2M}\eps \geq \frac{M}{M + \frac{1}{100}}.\]
It follows that $(f_k)_{\Z} = 0$ for all but at most $M + \frac{1}{100}$ values of $k$, as desired. \endproof

\section{Bourgain systems}

We now begin assembling the tools required to prove Lemma \ref{ind-step}.

Many theorems in additive combinatorics can be stated for an arbitrary abelian group $G$, but are much easier to prove in certain \emph{finite field models}, that is to say groups $G = \F_p^n$ where $p$ is a small fixed prime. This phenomenon is discussed in detail in the survey \cite{green-fin-field}. The basic reason for it is that the groups $\F_p^n$ have a very rich subgroup structure, whereas arbitrary groups need not: indeed the group $\Z/N\Z$, $N$ a prime, has no non-trivial subgroups at all. 

In his work on 3-term arithmetic progressions Bourgain \cite{bourgain} showed that \emph{Bohr sets} may be made to play the r\^ole of ``approximate subgroups'' in many arguments. A definition of Bohr sets will be given later. Since his work, similar ideas have been used in several places \cite{green-regularity,green-kony,green-tao-inverseu3,sanders1, sanders2,shkredov}.

In this paper we need a notion of approximate subgroup which includes that of Bourgain but is somewhat more general. In particular we need a notion which is invariant under Freiman isomorphism. A close examination of Bourgain's arguments reveals that the particular structure of Bohr sets is only relevant in one place, where it is necessary to classify the set of points at which the Fourier transform of a Bohr set is large. In an exposition of Bourgain's work, Tao \cite{tao:bourgain} showed how to do without this information, and as a result of this it is possible to proceed in more abstract terms.

\begin{definition}[Bourgain systems]\label{bour-def}
Let $G$ be a finite abelian group and let $d \geq 1$ be an integer. A Bourgain system $\S$ of dimension $d$ is a collection $(X_{\rho})_{\rho \in [0,4]}$ of subsets of $G$ indexed by the nonnegative real numbers such that the following axioms are satisfied:
\begin{enumerate}
\item[\textsc{bs1}] (Nesting) If $\rho' \leq \rho$ we have $X_{\rho'} \subseteq X_{\rho}$;
\item[\textsc{bs2}] (Zero) $0 \in X_0$;
\item[\textsc{bs3}] (Symmetry) If $x \in X_{\rho}$ then $-x \in X_{\rho}$;
\item[\textsc{bs4}] (Addition) For all $\rho,\rho'$ such that $\rho + \rho' \leq 4$ we have $X_{\rho} + X_{\rho'} \subseteq X_{\rho + \rho'}$;
\item[\textsc{bs5}] (Doubling) If $\rho \leq 1$ we have $|X_{2\rho}| \leq 2^d|X_{\rho}|$.
\end{enumerate}
We refer to $|X_1|$ as the \emph{size} of the system $\S$, and write $|\S|$ for this quantity. 
\end{definition}
\begin{remarks}
If a Bourgain system has dimension at most $d$, then it also has dimension at most $d'$ for any $d' \geq d$. It is convenient, however, to attach a fixed dimension to each system. Note that the definition is largely independent of the group $G$, a feature which enables one to think of the basic properties of Bourgain systems without paying much attention to the underlying group.
\end{remarks}

\begin{definition}[Measures on a Bourgain system]
Suppose that $\S = (X_{\rho})_{\rho \in [0,4]}$ is a Bourgain system contained in a group $G$. We associate to $\S$ a system $(\beta_{\rho})_{\rho \in [0,2]}$ of probability measures on the group $G$. These are defined by setting
\[ \beta_{\rho} := \frac{1_{X_{\rho}}}{|X_{\rho}|} \ast \frac{1_{X_{\rho}}}{|X_{\rho}|}.\]
Note that $\beta_{\rho}$ is supported on $X_{2\rho}$. 
\end{definition}

\begin{definition}[Density]
We define $\mu(\S) = |\S|/|G|$ to be the \textit{density} of $\S$ relative to $G$.
\end{definition}
 
\begin{remarks} 
Note that everything in these two definitions \emph{is} rather dependent on the underlying group $G$. The reason for defining our measures in this way is that the Fourier transform $\widehat{\beta}_{\rho}$ is real and nonnegative. This positivity property will be very useful to us later. The idea of achieving this by convolving an indicator function with itself goes back, of course, to Fej\'er. For a similar use of this device see \cite{green-regularity}, especially Lemma 7.2.
\end{remarks}

The first example of a Bourgain system is a rather trivial one.
\begin{example}[Subgroup systems]
Suppose that $H \leq G$ is a subgroup. Then the collection $(X_{\rho})_{\rho \in [0,4]}$ in which each $X_{\rho}$ is equal to $H$ is a Bourgain system of dimension $0$.
\end{example}

The second example is important only in the sense that later on it will help us economise on notation.
\begin{example}[Dilated systems]
Suppose that $\S = (X_{\rho})_{\rho \in [0, 4]}$ is a Bourgain system of dimension $d$. Then, for any $\lambda \in (0,1]$, so is the collection $\lambda \S := (X_{\lambda \rho})_{\rho \in [0,4]}$.
\end{example}

The following simple lemma concerning dilated Bourgain systems will be useful in the sequel.

\begin{lemma}\label{dilate-sizes}
Let $\S$ be a Bourgain system of dimension $d$, and suppose that $\lambda \in (0,1]$. Then $\dim(\lambda \S) = d$ and $|\lambda \S| \geq (\lambda/2)^d|\S|$.\endproof
\end{lemma}

\begin{definition}[Bohr systems]\label{bohr-def}  The first substantial example of a Bourgain system is the one contained in the original paper \cite{bourgain}. Let $\Gamma = \{\gamma_1,\dots,\gamma_k\} \subseteq \widehat{G}$ be a collection of characters, let $\kappa_1,\dots,\kappa_k > 0$, and define the system $\Bohr_{\kappa_1,\dots,\kappa_k}(\Gamma)$ by taking
\[ X_{\rho} := \{x \in G: |1 - \gamma_j(x)| \leq \kappa_j\rho\}.\]
When all the $\kappa_i$ are the same we write $\Bohr_{\kappa}(\Gamma) = \Bohr_{\kappa_1,\dots,\kappa_k}(\Gamma)$ for short.
The properties \textsc{bs1,bs2} and \textsc{bs3} are rather obvious. Property \textsc{bs4} is a consequence of the triangle inequality and the fact that $|\gamma(x) - \gamma(x')| = |1 - \gamma(x-x')|$. Property \textsc{bs5} and a lower bound on the density of Bohr systems are documented in the next lemma, a proof of which may be found in any of \cite{green-regularity,green-kony,green-tao-inverseu3}.
\end{definition}

\begin{lemma}\label{bohr-facts}
Suppose that $\S = \Bohr_{\kappa_1,\dots,\kappa_k}(\Gamma)$ is a Bohr system. Then $\dim(\S) \leq 3k$ and $|\S| \geq 8^{-k}\kappa_1\dots\kappa_k|G|$.\endproof
\end{lemma}

The notion of a Bourgain system is invariant under Freiman isomorphisms.

\begin{example}[Freiman isomorphs]
Suppose that $\S = (X_{\rho})_{\rho \in [0,4]}$ is a Bourgain system and that $\phi : X_4 \rightarrow G'$ is some Freiman isomorphism such that $\phi(0) = 0$. Then $\phi(\S) := (\phi(X_{\rho}))_{\rho \in [0,4]}$ is a Bourgain system of the same dimension and size.
\end{example}

The next example is of no real importance over and above those already given, but it does serve to set the definition of Bourgain system in a somewhat different light.

\begin{example}[Translation-invariant pseudometrics]
Suppose that $d : G \times G \rightarrow \R_{\geq 0}$ is a translation-invariant pseudometric. That is, $d$ satisfies the usual axioms of a metric space except that it is possible for $d(x,y)$ to equal zero when $x \neq y$ and we insist that $d(x + z,y+z) = d(x,y)$ for any $x,y,z$. Write $X_{\rho}$ for the ball
\[ X_{\rho} := \{x \in G : d(x,0) \leq \rho\}.\]
Then $(X_{\rho})_{\rho \in [0,4]}$ is a Bourgain system precisely if $d$ is \emph{doubling}, cf. \cite[Ch. 1]{heinonen}.
\end{example}

\begin{remark}
It might seem more elegant to try and define a Bourgain system to be the same thing as a doubling, translation invariant pseudometric. There is a slight issue, however, which is that such Bourgain systems satisfy \textsc{bs1}--\textsc{bs5} for all $\rho \in [0,\infty)$. It is not in general possible to extend a Bourgain system defined for $\rho \in [0,4]$ to one defined for all nonnegative $\rho$, as one cannot keep control of the dimension condition \textsc{bs5}. Consider for example the (rather trivial) Bourgain system in which $X_{\rho} = \{0\}$ for $\rho < 4$ and $X_4$ is a symmetric set of $K$ ``dissociated'' points. 
\end{remark}

We now proceed to develop the basic theory of Bourgain systems. For the most part this parallels the theory of Bohr sets as given in several of the papers cited earlier. The following lemmata all concern a Bourgain system $\S$ with dimension $d$.

We begin with simple covering and metric entropy estimates. The following covering lemma could easily be generalized somewhat, but we give here just the result we shall need later on.

\begin{lemma}[Covering lemma]\label{covering}
For any $\rho \leq 1/2$, $X_{2\rho}$ may be covered by $2^{4d}$ translates of $X_{\rho/2}$.
\end{lemma} \proof Let $Y = \{y_1,\dots,y_k\}$ be a maximal collection of elements of $X_{2\rho}$ with the property that the balls $y_j + X_{\rho/4}$ are all disjoint. If there is some point $y \in X_{2\rho}$ which does not lie in any $y_j + X_{\rho/2}$, then $y + X_{\rho/4}$ does not intersect $y_j + X_{\rho/4}$ for any $j$ by \textsc{bs4}, contrary to the supposed maximality of $Y$. Now another application of \textsc{bs4} implies that
\[ \bigcup_{j = 1}^k (y_j + X_{\rho/4}) \subseteq X_{9\rho/4}.\]
We therefore have 
\[ k \leq |X_{9\rho/4}|/|X_{\rho/4}| \leq |X_{4\rho}|/|X_{\rho/4}| \leq 2^{4d}.\]
The lemma follows.\endproof

\begin{lemma}[Metric entropy lemma]\label{entropy}
Let $\rho \leq 1$. The group $G$ may be covered by at most $(4/\rho)^d \mu( \S )^{-1}$ translates of $X_{\rho}$. 
\end{lemma} \proof This is a simple application fo the Ruzsa covering lemma (cf. \cite[Ch. 2]{tao-vu-book}) and the basic properties of Bourgain systems. Indeed the Ruzsa covering lemma provides a set $T \subseteq G$ such that $G = X_{\rho/2} - X_{\rho/2} + T$, where
\[ |T| \leq \frac{|X_{\rho/2} + G|}{|X_{\rho/2}|} \leq \frac{|G|}{|X_{\rho/2}|}.\] \textsc{bs4} then tells us that $G = X_{\rho} + T$. To bound the size of $T$ above, we observe from \textsc{bs5} that $|X_{\rho/2}| \geq (\rho/4)^d |X_1|$. The result follows.\endproof

In this paper we will often be doing a kind of Fourier analysis relative to Bourgain systems. In this regard it is useful to know what happens when an arbitrary Bourgain system $(X_{\rho})_{\rho \in [0,4]}$ is joined with a system $(\Bohr(\Gamma,\eps\rho))_{\rho \in [0,4]}$ of Bohr sets, where $\Gamma \subseteq \widehat{G}$ is a set of characters. It turns out not to be much harder to deal with the join of a pair of Bourgain systems in general. 

\begin{definition}[Joining of two Bourgain systems]
Suppose that $\S = (X_{\rho})_{\rho \in [0,4]}$ and $\S' = (X'_{\rho})_{\rho \in [0,4]}$ are two Bourgain systems with dimensions at most $d$ and $d'$ respectively. Then we define the join of $\S$ and $\S'$, $\S \wedge \S'$, to be the collection $(X_{\rho} \cap X'_{\rho})_{\rho \in [0,4]}$.
\end{definition}

\begin{lemma}[Properties of joins]\label{join-lem}
Let $\S,\S'$ be as above. Then the join $\S \wedge \S'$ is also a Bourgain system. It has dimension at most $4(d+d')$ and its size satisfies the bound
\[ |\S \wedge \S'| \geq 2^{-3(d + d')}\mu(\S') |\S |.\]
\end{lemma}
\proof It is trivial to verify properties \textsc{bs1}--\textsc{bs4}. To prove \textsc{bs5}, we apply Lemma \ref{covering} to both $\S$ and $\S'$. This enables us to cover $X_{2\rho} \cap X'_{2 \rho}$ by at $2^{4(d + d')}$ sets of the form $T = (y + X_{\rho/2}) \cap (y' + X'_{ \rho/2})$. Now for any fixed $t_0 \in T$ the map $t \mapsto t - t_0$ is an injection from $T$ to $X_{\rho} \cap X'_{\rho}$. It follows, of course, that $|T| \leq |X_{\rho} \cap X'_{\rho}|$ and hence that
\[ |X_{2\rho} \cap X'_{2\rho}| \leq 2^{4(d + d')}|X_{\rho} \cap X'_{\rho}|.\]
This establishes the claimed bound on the dimension of $\S \wedge \S'$. It remains to obtain a lower bound for the density of this system. To do this, we apply Lemma \ref{entropy} to cover $G$ by at most $8^{d'}\mu(\S')^{-1}$ translates of $X'_{1/2}$. It follows that there is some $x$ such that \[ |X_{1/2} \cap (x + X'_{1/2})| \geq 8^{-d'}\mu(\S')^{-1}|X_{1/2}| \geq 2^{-3(d+d')}\mu(\S') |X_1|.\] Now for any fixed $x_0 \in X_{1/2} \cap (x + X'_{1/2})$ the map $x \mapsto x - x_0$ is an injection from $X_{1/2} \cap (x + X'_{1/2})$ to $X_1 \cap X'_{1}$. It follows that
\[ |X_1 \cap X'_{1}| \geq 2^{-3(d + d')}\mu(\S') |X_1|,\] which is equivalent to the lower bound on the size of $\S \wedge \S'$ that we claimed.\endproof

We move on now to one of the more technical aspects of the theory of Bourgain systems, the notion of regularity.

\begin{definition}[Regular Bourgain systems]
Let $\S = (X_{\rho})_{\rho \in [0,4]}$ be a Bourgain system of dimension $d$. We say that the system is \emph{regular} if
\[ 1 - 10 d\kappa \leq \frac{|X_{1}|}{|X_{1 + \kappa}|} \leq 1 + 10 d \kappa\] whenever $|\kappa| \leq 1/10d$.
\end{definition}

\begin{lemma}[Finding regular Bourgain systems]\label{ubiq-lem}
Suppose $\S$ is a Bourgain system. Then there is some $\lambda \in [1/2,1]$ such that the dilated system $\lambda \S$ is regular.
\end{lemma}
\proof Let $f : [0,1] \rightarrow \R$ be the function $f(a) := \frac{1}{d}\log_2 |X_{2^a}|$. Observe that $f$ is non-decreasing in $a$ and that $f(1) - f(0) \leq 1$. We claim that there is an $a \in [\frac{1}{6}, \frac{5}{6}]$ such that $|f(a+ x) - f(a)| \leq 3|x|$ for all $|x| \leq \frac{1}{6}$. If no such $a$ exists then for every $a \in [\frac{1}{6}, \frac{5}{6}]$ there is an interval $I(a)$ of length at most $\frac{1}{6}$ having one endpoint equal to $a$ and with $\int_{I(a)} df > \int_{I(a)} 3 dx$. These intervals cover $[\frac{1}{6}, \frac{5}{6}]$, which has total length $\frac{2}{3}$. A simple covering lemma that has been discussed by Hallard Croft \cite{croft-eureka} (see also \cite[Lemma 3.4]{green-kony}) then allows us to pass to a disjoint subcollection $I_1 \cup \dots \cup I_n$ of these intervals with total length at least $\frac{1}{3}$. However we now have
\[ 1 \geq \int^1_0 df \geq \sum_{i=1}^n \int_{I_i} df > \sum_{i = 1}^n \int_{I_i} 3 \, dx \geq \frac{1}{3} \cdot 3, \] a contradiction.
It follows that there is indeed an $a$ such that $|f(a+x) - f(a)| \leq 3|x|$ for all $|x| \leq \frac{1}{6}$. Setting $\lambda := 2^a$, it is easy to see that
\[ e^{-5d  \kappa} \leq \frac{|X_{\lambda}|}{|X_{(1 + \kappa)\lambda}|} \leq e^{5d \kappa}\] whenever $|\kappa| \leq 1/10d$. Since $1 - 2|x| \leq e^x \leq 1 + 2|x|$ when $|x| \leq 1$, it follows that $\lambda \S$ is a regular Bourgain system.\endproof

\begin{lemma}\label{lemma3.2}
Suppose that $\S$ is a regular Bourgain system of dimension $d$ and let $\kappa \in (0,1)$. Suppose that $y \in X_{\kappa}$. Then
\[ \E_{x \in G} |\beta_1(x + y) - \beta_1(x)| \leq 20d \kappa.\] 
\end{lemma} \proof For this lemma only, let us write $\mu_1 := 1_{X_1}/|X_1|$, so that $\beta_1 = \mu_1 \ast \mu_1$. We first claim that if $y \in X_{\kappa}$ then
\[ \E_{x \in G} |\mu_1(x + y) - \mu_1(x)| \leq 20 d \kappa.\]
The result is trivial if $\kappa > 1/10d$, so assume that $\kappa \leq 1/10d$. Observe that $|\mu_1(x + y) - \mu_1(x)| = 0$ unless $x \in X_{1 + \kappa} \setminus X_{1 - \kappa}$. Since $\S$ is regular, the size of this set is at most $20 d \kappa |X_{1}|$, and the claim follows immediately.

To prove the lemma, note that
\begin{align*}
\E_{x \in G} |\beta_1(x+y) - \beta_1(x)| &= \E_{x \in G} |\mu_1 \ast \mu_1 (x+y) - \mu_1\ast \mu_1(x)| \\ & = \E_x \big| \E_z \mu_1(z) \mu_1(x + y - z) - \E_z \mu_1(z) \mu_1(x - z)\big| \\ & \leq  \E_z \mu_1(z) \E_x | \mu_1(x+y - z) - \mu_1(x - z)| \\ & \leq 20d\kappa,
\end{align*}
the last inequality following from the claim.\endproof

The operation of convolution by $\beta_1$ will play an important r\^ole in this paper, particularly in the next section. 

\begin{definition}[Convolution operator]
Suppose that $\S$ is a Bourgain system. Then we associate to $\S$ the map $\psi_{\S} : L^{\infty}(G) \rightarrow L^{\infty}(G)$ defined by $\psi_{\S}f := f \ast \beta_1$, or equivalently by $(\psi_{\S}f)^{\wedge} := \widehat{f} \widehat{\beta}_1$. 
\end{definition}

We note in particular that, since $\widehat{\beta}_1$ is real and non-negative, we have
\begin{equation}\label{fund-property} \Vert f \Vert_A = \Vert \psi_{\S}f \Vert_A + \Vert f - \psi_{\S}f \Vert_A.\end{equation}

\begin{lemma}[Almost invariance]\label{lemma3.3}
Let $f : G \rightarrow \C$ be any function. Let $\S$ be a regular Bourgain system of dimension $d$, let $\kappa \in (0,1)$ and suppose that $y \in X_{\kappa}$. Then \[ |\psi_{\S}f(x + y) - \psi_{\S}f(x)| \leq 20 d \kappa \Vert f \Vert_{\infty}\] for all $x \in G$. 
\end{lemma}
\proof The left hand side, written out in full, is
\[| \E_t f(t) (\beta_{\rho}(t - x - y) - \beta_{\rho}(t - x))|.\]
The lemma follows immediately from Lemma \ref{lemma3.2} and the triangle inequality.\endproof

\begin{lemma}[Structure of Spec]\label{spec-struct} Let $\delta \in (0,1]$.
Suppose that $\S$ is a regular Bourgain system of dimension $d$ and that $\gamma \in \Spec_{\delta}(\beta_{1})$. Suppose that $\kappa \in (0,1)$. Then we have 
\[ |1 - \gamma(y)| \leq 20\kappa d/\delta\]
for all $y \in X_{\kappa}$.
\end{lemma}
\proof Suppose that $y \in X_{\kappa}$. Then we have
\begin{align*} \delta |1 - \gamma(y)| \leq |\widehat{\beta}_{1}(\gamma)| |1 - \gamma(y)|& = | \E_{x \in G} \beta_{1}(x) (\gamma(x) - \gamma(x + y))| \\ & = |\E_{x \in G} (\beta_{1}(x) - \beta_{1}(x - y)) \gamma(x)|.\end{align*}
This is bounded by $20 d \kappa$ by Lemma \ref{lemma3.2}.\endproof

\section{Averaging over a Bourgain system}

Let $\S = (X_{\rho})_{\rho \in [0,4]}$ be a Bourgain system of dimension $d$. Recall from the last section the definition of the operator $\psi_{\S} : L^{\infty}(G) \rightarrow L^{\infty}(G)$.
Recalling our earlier paper \cite{green-sanders}, one might use operators of this type to effect a decomposition
\[ f = \psi_{\S}f + (f - \psi_{\S}f),\]
the aim being to prove Theorem \ref{mainthm-finite} by induction on $\Vert f \Vert_A$. To make such a strategy work, a judicious choice of $\S$ must be made. First of all it must be ensured that both $\Vert\psi_{\S} f\Vert_A$ and $\Vert f - \psi_{\S}f \Vert_A$ are significantly less than $\Vert f \Vert_A$. In this regard \eqref{fund-property} is of some importance, and this is why we defined the measures $\beta_{\rho}$ in such a way that $\widehat{\beta}_{\rho}$ is always real and nonnegative. The actual accomplishment of this will be a task for the next section. In an ideal world, our second requirement would be that $\psi_{\S}$ preserves the property of being $\Z$-valued. As in our earlier paper this turns out not to be possible and one must expand the collection of functions under consideration to include those for which $d(f,\Z) \leq \eps$. The reader may care to recall the definition of $d(f,\Z)$ and of $f_{\Z}$ at this point: they are given at the start of Section \ref{main-deduction-sec}.

The main result of this section states that if $f$ is almost integer-valued then any Bourgain system $\S$ may be refined to a system $\S'$ so that $\psi_{\S'}f$ is almost integer-valued. A result of this type in the finite field setting, where $\S$ is just a subgroup system in $\F_2^n$, was obtained in \cite{green-sanders}. The argument there, which was a combination of \cite[Lemma 3.4]{green-sanders} and \cite[Proposition 3.7]{green-sanders}, was somewhat elaborate and involved polynomials which are small near small integers. The argument we give here is different and is close to the main argument in \cite{green-kony} (in fact, it is very close to the somewhat simpler argument, leading to a bound of $O(\log^{-1/4} p)$, sketched just after Lemma 4.1 of that paper). In the finite field setting it is simpler than that given in \cite[Sec. 3]{green-sanders} and provides a better bound. Due to losses elsewhere in the argument, however, it does not lead to an improvement in the overall bound in our earlier paper.

\begin{proposition}\label{prop4.2}
Suppose that $f : G \rightarrow \R$ satisfies $\Vert f \Vert_A \leq M$, where $M \geq 1$, and also $d(f,\Z) < 1/4$. Let $\S$ be a regular Bourgain system of dimension $d \geq 2$, and let $\eps \leq \frac{1}{4}$ be a positive real. Then there is a regular Bourgain system $\S'$ with dimension $d'$ such that
\begin{equation}\label{dim-bound}
d' \leq  4d + \frac{64M^2}{\eps^2};
\end{equation}
\begin{equation}\label{size-bd}| \S' | \geq e^{-\frac{CdM^4}{\eps^4}\log(dM/\eps)}| \S|;\end{equation} 
\begin{equation}\label{avg-lwr} \Vert \psi_{\S'}f \Vert_{\infty} \geq \Vert \psi_{\S} f \Vert_{\infty} - \eps\end{equation}
and such that
\begin{equation}\label{almost-int-bd} d(\psi_{\S'}f,\Z) \leq d(f,\Z) + \eps.\end{equation}
\end{proposition}
\begin{remarks}
The stipulation that $d \geq 2$ and that $M \geq 1$ is made for notational convenience in our bounds. These conditions may clearly be satisfied in any case by simply increasing $d$ or $M$ as necessary.
\end{remarks}
\proof We shall actually find $\S'$ satisfying the following property:
\begin{equation}\label{avg-bd} \E_{x \in G} (f - \psi_{\S'}f)(x)^2 \beta'_{\rho}(x- x_0) \leq \eps^2/4\end{equation} for any $x_0 \in G$ and any every $\rho \geq \eps/160d'M$ such that $\rho \S'$ is regular. The truth of \eqref{avg-bd} implies \eqref{almost-int-bd}. To see this, suppose that \eqref{almost-int-bd} is false. Then there is $x_0$ so that $\psi_{\S'}f(x_0)$ is not within $d(f,\Z) + \eps$ of an integer. Noting that $\Vert f \Vert_{\infty} \leq \Vert f \Vert_A \leq M$, we see from Lemma \ref{lemma3.3} that $\psi_{\S'}f(x)$ is not within $d(f,\Z) + \eps/2$ of an integer for any $x \in x_0 + X_{\eps/40d'M}$. Choosing, according to Lemma \ref{ubiq-lem}, a value $\rho \in [\eps/160d'M,\eps/80d'M]$ for which $\rho \S'$ is regular, we have
\[ \E_x (f - \psi_{\S'}f)(x)^2 \beta'_{\rho}(x - x_0) > \eps^2/4,\] contrary to our assumption of \eqref{avg-bd}.

It remains to find an $\S'$ such that \eqref{avg-bd} is satisfied for all $x_0 \in G$ and all $\rho \geq \eps/160d'M$ such that $\rho \S'$ is regular.
We shall define a nested sequence $\S^{(j)} = (X_{\rho}^{(j)})_{\rho \in [0,4]}$, $j = 0,1,2,\dots$ of regular Bourgain systems with $d_j := \dim(\S^{(j)})$. We initialize this process by taking $\S^{(0)} := \S$.

Suppose that $\S^{(j)}$ does not satisfy \eqref{avg-bd}, that is to say there is $y \in G$ and $\rho \geq \eps/160d_jM$ such that $\rho \S^{(j)}$ is regular and
\[ \E_{x \in G} (f - f\ast \beta_1^{(j)})(x)^2 \beta^{(j)}_{\rho}(x - y) > \eps^2/4.\]
Applying Plancherel, this means that
\[ \sum_{\gamma \in \widehat{G}} \big( (f - f \ast \beta_1^{(j)})\beta^{(j)}_{\rho}(\cdot - y)\big)^{\wedge}(\gamma) (f - f\ast \beta_1^{(j)})^{\wedge}(\gamma) > \eps^2/4.\] This implies that
\[ \Vert \big( (f - f\ast \beta_1^{(j)})\beta^{(j)}_{\rho}( \cdot - y)\big)^{\wedge} \Vert_{\infty} > \eps^2/8M,\]
which implies as that there is some $\gamma_0^{(j+1)} \in \widehat{G}$ such that
\[  \sum_{\gamma} |\widehat{f}(\gamma)| |1 - \widehat{\beta}_1^{(j)}(\gamma)| |\widehat{\beta}^{(j)}_{\rho} (\gamma_0^{(j+1)} - \gamma)| > \eps^2/8M.\] Removing the tails where either $|1 - \widehat{\beta}_1^{(j)}(\gamma)| \leq \eps^2/32M^2$ or $|\widehat{\beta}^{(j)}_{\rho}(\gamma_0^{(j+1)} - \gamma)| \leq \eps^2/64M^2$, this implies that
\begin{equation}\label{l1-mass} \sum_{\gamma \in \Gamma^{(j)}} |\widehat{f}(\gamma)| > \eps^2/16M,\end{equation}
where the sum is over the set \[\Gamma^{(j)} := \big(\gamma_0^{(j+1)} + \Spec_{\eps^2/64M^2}(\beta^{(j)}_{\rho})\big)\setminus \Spec_{1 - \eps^2/32M^2}(\beta^{(j)}_1).\]
We shall choose a regular Bourgain system $\S^{(j+1)}$ in such a way that
\begin{equation}\label{to-sat}
\gamma_0^{(j+1)} + \Spec_{\eps^2/64M^2}(\beta^{(j)}_{\rho}) \subseteq \Spec_{1 - \eps^2/32M^2}(\beta^{(j+1)}_1).
\end{equation}
The sets $\Gamma^{(j)}$ are then disjoint, and it follows from \eqref{l1-mass} that the iteration must stop for some $j = J$, $J \leq 16M^2/\eps^2$. We then define $\S' := \S^{(J)}$.

To satisfy \eqref{to-sat} we take
\begin{equation}\label{sj-def} \S^{(j+1)} := \lambda\big(\kappa \rho \S^{(j)} \wedge \Bohr_{\kappa'}(\{\gamma_0^{(j+1)}\})\big),\end{equation}
where $\kappa := 2^{-17}\epsilon^4/d_jM^4$, $\kappa' := \eps^2/64 M^2$, and $\lambda \in [1/2,1]$ is chosen so that $\S^{(j+1)}$ is regular. Note that
\begin{equation}\label{eq800}
\lambda \kappa \rho \geq \frac{\eps^5}{2^{26}d_j^2 M^5}.
\end{equation}
Suppose that $\gamma \in \Spec_{\eps^2/64M^2}(\beta_{\rho}^{(j)})$. Then in view of Lemma \ref{spec-struct} and the fact that $\rho \S^{(j)}$ is regular we have
\[ |1 - \gamma(x)| \leq \frac{1280 \kappa d_j M^2}{\eps^2} \leq \frac{\eps^2}{64 M^2}\] whenever $x \in X_1^{(j+1)}$. Furthermore we also have
\[ |1 - \gamma_0^{(j+1)}(x)| \leq \kappa' = \eps^2/64M^2.\]
It follows that if $x \in X_1^{(j+1)}$ then
\[ |1 - \gamma_0^{(j+1)}(x)\gamma(x)| \leq \eps^2/32M^2,\]
and therefore $\gamma_0^{(j+1)} + \gamma \in \Spec_{1 - \eps^2/32M^2}(\beta_1^{(j+1)})$.
 
It remains to bound $\dim(\S^{(j)})$ and $|\S^{(j)}|$. To this end we note that by construction we have
\[ \S^{(j)} = \delta^{(j)} \S^{(0)} \wedge \Bohr_{\kappa_1,\dots,\kappa_j}(\{\gamma_0^{(1)},\dots,\gamma_0^{(j)}\}),\] where each $\kappa_i$ is at least $2^{-j}\epsilon^2/64M^2$ and, in view of \eqref{eq800},
\[ \delta^{(j)} \geq (\eps^5/2^{26}M^5)^j \big(\prod_{i \leq j} d_i\big)^{-2}.\]
It follows from Lemmas \ref{bohr-facts} and \ref{join-lem} that $d_j \leq 4(d + j)$ for all $j$, and in particular we obtain the claimed upper bound on $\dim(\S')$. It follows from the same two lemmas together with Lemma \ref{dilate-sizes} and a short computation that $| \S' |$ is subject to the claimed lower bound. The lower bound we have given is, in fact, rather crude but has been favoured due to its simplicity of form.

It remains to establish \eqref{avg-lwr}. 
Noting that 
\[ f \ast \beta_1  = f \ast \beta_1\ast \beta'_1 - f \ast ( \beta_1 \ast \beta'_1 - \beta_1),\]
we obtain the bound
\[ \Vert \psi_{\S} f \Vert_{\infty}  \leq \Vert \psi_{\S'}f\Vert_{\infty} + \Vert f \Vert_{\infty} \Vert \beta_1 \ast \beta'_1 - \beta_1 \Vert_{1}.\]
If 
\[ \Vert \beta_1 \ast \beta'_1 - \beta_1 \Vert_{1} \leq \eps/M\]
then the result will follow. We have, however, that
\[ \Vert \beta_1 \ast \beta'_1 - \beta_1 \Vert_1 = \E_x \big|\E_y \beta'_1(y) (\beta_1(x-y) - \beta_1(x))  \big|,\] and from Lemma \ref{lemma3.2} it will follow that this is at most $\eps/M$ provided that $\Supp(\beta'_1) \subseteq X_{\eps/20dM}$. This, however, is more than guaranteed by the construction of the successive Bourgain systems as given in \eqref{sj-def}. Note that we may assume without loss of generality that the iteration does proceed for at least one step: even if \eqref{avg-bd} is satisfied by $\S = \S^{(0)}$, we may simply take an arbitrary $\gamma_0^{(1)} \in \widehat{G}$ and define $\S^{(1)}$ as in \eqref{sj-def}.
\endproof

\section{A weak Freiman theorem}

In our earlier work \cite{green-sanders} we used (a refinement of) Ruzsa's analogue of Freiman's theorem, which gives a fairly strong characterisation of subsets $A \subseteq \F_2^n$ satisfying a small doubling condition $|A + A| \leq K|A|$. An analogue of this theorem for any abelian group was obtained in \cite{green-ruzsa}. We could apply this theorem here, but as reward for setting up the notion of Bourgain systems in some generality we are able to make do with a weaker theorem of the following type, which we refer to as a ``weak Freiman theorem''.

\begin{proposition}[Weak Freiman]\label{weak-frei}
Suppose that $G$ is a finite abelian group, and that $A \subseteq G$ is a finite set with $|A + A| \leq K|A|$. Then there is a regular Bourgain system $\S = (X_{\rho})_{\rho \in [0,4]}$ such that
\[ \dim(\S) \leq CK^C;\]
\[ |\S|  \geq e^{-CK^C} |A|\] and
\[ \Vert \psi_{\S} 1_A \Vert_{\infty} \geq cK^{-C}.\]
\end{proposition}

\begin{remark}
We note that, unlike the usual Freiman theorem, it is clear how one might formulate a weak Freiman theorem in arbitrary (non-abelian) groups. We are not able to prove such a statement, and there seem to be significant difficulties in doing so. For example, there is no analogue of \cite[Proposition 1.2]{green-ruzsa} in general groups. See \cite{green-nonabelian} for more details.
\end{remark}

We begin by proving a result similar to Proposition \ref{weak-frei} in what appears to be a special case: when $A$ is a dense subset of a group $G$. We will show later on that the general case can be reduced to this one.

\begin{proposition}[Bogolyubov-Chang argument]\label{bog-chang} Let $G$ be a finite abelian group, and suppose that $A \subseteq G$ be a set with $|A| = \alpha |G|$ and $|A + A| \leq K|A|$. Then there is a regular Bourgain system $\S$ with
\[ \dim(\S) \leq CK\log(1/\alpha),\]
\[ \Vert \psi_{\S} 1_A \Vert_{\infty} \geq 1/2K,\]
\[ |\S| \geq (CK\log(1/\alpha))^{-CK\log(1/\alpha)}|G|\]
and 
\[ X_4 \subseteq 2A - 2A.\]
\end{proposition}
\proof The argument is a variant due to Chang \cite{chang-freiman} of an old argument of Bogolyubov \cite{bog}. It is by now described in several places, including the book \cite{tao-vu-book}.

Set 
\[ \Gamma := \Spec_{1/4\sqrt{K}}(A) := \{\gamma \in \widehat{G}: |\widehat{1}_A(\gamma)| \geq \frac{\alpha}{4\sqrt{K}}\},\]
and take 
\[ \tilde{\S} = (\tilde{X}_{\rho})_{\rho \in [0,4]} := \Bohr_{1/20}(\Gamma),\] a Bohr system as defined in Definition \ref{bohr-def}. We claim that $\tilde{X}_4 \subseteq 2A - 2A$. Recall from the definition that $\tilde{X}_4$ consists of those $x \in G$ for which $|1 - \gamma(x)| \leq \frac{1}{5}$ for all $\gamma \in \Gamma$. Suppose then that $x \in \tilde{X}_4$. By the inversion formula we have
\begin{align*}
\Vert \widehat{1}_A \Vert_4^4 - 1_A \ast 1_A \ast 1_{-A} \ast 1_{-A}(x) &= \sum_{\gamma} |\widehat{1}_A(\gamma)|^4 (1 - \gamma(x)) \\ & \leq \sum_{\gamma \in \Gamma} |\widehat{1}_A(\gamma)|^4 |1 - \gamma(x)| + \sum_{\gamma \notin \Gamma} |\widehat{1}_A(\gamma)|^4 |1 - \gamma(x)|\\
 & \leq \frac{1}{5} \Vert \widehat{1}_A \Vert_4^4 + \frac{\alpha^2}{8K} \Vert 1_A \Vert_2^2 \\ & = \frac{1}{5} \Vert \widehat{1}_A \Vert_4^4 + \frac{\alpha^3}{8K}.
\end{align*}
However the fact that $|A + A| \leq K|A|$ implies, using Cauchy-Schwarz, that
\begin{equation}\label{cbs-app} \Vert \widehat{1}_A \Vert_4^4 = \Vert 1_A \ast 1_A \Vert_2^2 \geq \alpha^3/K.\end{equation} It follows that
\[ \Vert \widehat{1}_A \Vert_4^4 - 1_A \ast 1_A \ast 1_{-A} \ast 1_{-A}(x) \leq (\frac{1}{5} + \frac{1}{8}) \Vert \widehat{1}_A \Vert_4^4 < \Vert \widehat{1}_A \Vert_4^4.\]
Therefore $1_A \ast 1_A \ast 1_{-A} \ast 1_{-A}(x) > 0$, that is to say $x \in 2A - 2A$.

Now we only have the dimension bound $\dim(\tilde{\S}) \leq 48K/\alpha$, which comes from the fact (a consequence of Parseval's identity) that $|\Gamma| \leq 16K/\alpha$. This is substantially weaker than the bound $CK\log(1/\alpha)$ that we require. To obtain the superior bound we must refine $\tilde{\S}$ to a somewhat smaller system $\S$. To do this we apply a well-known lemma of Chang \cite{chang-freiman}, which follows from an inequality of Rudin \cite{rudin-book}. See also \cite{green-ruzsa,tao-vu-book} for complete, self-contained proofs of this result. In our case the lemma states that there is a set $\Lambda \subseteq \widehat{G}$, $|\Lambda| \leq 32K\log(1/\alpha)$, such that $\Gamma \subseteq \langle \Lambda \rangle$. Here,
\[ \langle \Lambda \rangle := \{\lambda_1^{\eps_1} \dots \lambda_k^{\eps_k} : \eps_i \in \{-1,0,1\}\},\]
where $\lambda_1,\dots,\lambda_k$ is a list of the characters in $\Lambda$. 

Now by repeated applications of the triangle inequality we see that
\[ \Bohr_{1/20k}(\Lambda) \subseteq \Bohr_{1/20}(\langle \Lambda \rangle) \subseteq \Bohr_{1/20}(\Gamma).\]
Thus if we set
\[ \S = (X_{\rho})_{\rho \in [0,4]} := \Bohr_{\lambda/20k}(\Lambda),\] where $\lambda \in [1/2,1]$ is chosen so that $\S$ is regular, then we have $X_4 \subseteq \tilde{X}_4 \subseteq 2A - 2A$. It follows from Lemma \ref{bohr-def} that $\dim(\S) \leq 72K\log (1/\alpha)$ and that
\[ |\S| \geq (1/320k)^k |G|  \geq (CK\log(1/\alpha))^{-CK\log(1/\alpha)}|G|,\] as required.

It remains to show that $\Vert \psi_{\S}1_A \Vert_{\infty} \geq 1/2K$. Let us begin by noting that if $\gamma \in \Gamma$ and $x \in \widetilde{X}_2$ then $|1 - \gamma(x)| \leq \frac{1}{10}$, and so if $\gamma \in \Gamma$ then $|\widehat{\beta}_1(\gamma)| \geq \frac{9}{10}$. It follows that
\begin{align*}
\Vert 1_A \ast 1_A \ast \beta_1 \Vert_2^2 &= \E 1_A \ast 1_A \ast \beta_1(x)^2 \\ &= \sum_{\gamma} |\widehat{1}_A(\gamma)|^4 |\widehat{\beta}_1(\gamma)|^2 \\ &\geq \sum_{\gamma \in \Gamma} |\widehat{1}_A(\gamma)|^4 |\widehat{\beta}_1(\gamma)|^2 - \frac{\alpha^3}{16K}\\ &\geq \frac{3}{4} \sum_{\gamma \in \Gamma} |\widehat{1}_A(\gamma)|^4 - \frac{\alpha^3}{16K} \\ & \geq \frac{3}{4}\big(\Vert \widehat{1}_A \Vert_4^4 - \frac{\alpha^3}{16K} \big) - \frac{\alpha^3}{16 K} \\ & \geq \alpha^3/2K,
\end{align*}
the last step following from \eqref{cbs-app}.
Since $\Vert 1_A \ast 1_A \ast \beta_1 \Vert_1 = \alpha^2$, it follows that 
\[ \Vert 1_A \ast 1_A \ast \beta_1 \Vert_{\infty} \geq \alpha/2K,\]
and hence that 
\[ \Vert \psi_{\S} 1_A \Vert_{\infty} = \Vert 1_A \ast \beta_1 \Vert_{\infty} \geq 1/2K,\]
as required.\endproof

\emph{Proof of Theorem \ref{weak-frei}}. By \cite[Proposition 1.2]{green-ruzsa} there is an abelian group $G'$, $|G'| \leq (CK)^{CK^2}|A|$, and a subset $A' \subseteq G'$ such that $A' \cong_{14} A$. We apply Proposition \ref{bog-chang} to this set $A'$. Noting that $\alpha \geq (CK)^{-CK^2}$, we obtain a Bourgain system $\S' = (X'_{\rho})_{\rho \in [0,4]}$ for which
\[ \dim(\S') \leq CK^C;\]
\[ |\S'| \geq e^{-CK^C} |A'|;\]
\[ \Vert \psi_{\S'} 1_{A'} \Vert_{\infty} \geq cK^{-C}\] and
\[ X'_4 \subseteq 2A' - 2A'.\]
Write $\phi : A' \rightarrow A$ for the Freiman 14-isomorphism between $A'$ and $A$. The map $\phi$ extends to a well-defined 1-1 map on $kA' - lA'$ for any $k,l$ with $k + l \leq 14$. By abuse of notation we write $\phi$ for any such map. In particular $\phi(0)$ is well-defined and we may define a ``centred'' Freiman 14-isomorphism $\phi_0(x) := \phi(x) - \phi(0)$.

Define $\S := \phi_0(\S')$. Since $X'_4 \subseteq 2A' - 2A'$, $\phi_0$ is a Freiman $2$-isomorphism on $X'_4$ with $\phi_0(0) = 0$. Therefore $\S$ is indeed a Bourgain system, with the same dimension and size as $\S'$.

It remains to check that $\Vert \psi_{\S} 1_A \Vert_{\infty} \geq cK^{-C}$. The fact that $\Vert \psi_{\S'} 1_{A'} \Vert_{\infty} \geq cK^{-C}$ means that there is $x$ such that $|1_{A'} \ast \beta'_1(x)| \geq cK^{-C}$. Since $\Supp(\beta'_1) \subseteq X'_2 \subseteq X'_4 \subseteq 2A' - 2A'$, we must have $x \in 3A' - 2A'$. We claim that $1_{A} \ast \beta_1(\phi(x)) = 1_{A'} \ast \beta'_1(x)$, which clearly suffices to prove the result. Recalling the definition of $\beta_1,\beta'_1$, we see that this amounts to showing that the number of solutions to
\[ x = a' - t'_1 + t'_2, \;\; a' \in A', t_i' \in X'_1,\]
is the same as the number of solutions to
\[ \phi_0(x) = \phi_0(a') - \phi_0(t'_1) + \phi_0(t'_2), \;\;a' \in A', t_i' \in X'_1.\] 
All we need check is that if $y \in 7A' - 7A'$ then $\phi_0(y) = 0$ only if $y = 0$. But since $0 \in 7A' - 7A'$, this follows from the fact that $\phi_0$ is 1-1 on $7A' - 7A'$.

To conclude the section we note that Proposition \ref{weak-frei} may be strengthened by combining it with the Balog-Szemer\'edi-Gowers theorem \cite[Proposition 12]{gowers-4aps} to obtain the following result.

\begin{proposition}[Weak Balog-Szemer\'edi-Gowers-Freiman]\label{weak-bsg}
Let $A$ be a subset of an abelian group $G$, and suppose that there are at least $\delta|A|^3$ additive quadruples $(a_1,a_2,a_3,$ $a_4)$ in $A^4$ with $a_1 + a_2 = a_3 + a_4$. Then there is a regular Bourgain system $\S$ satisfying
\[ \dim(\S) \leq C\delta^{-C};\]
\[ |\S| \geq e^{-C\delta^{-C}}|A|\]
and
\[ \Vert \psi_{\S}1_A \Vert_{\infty} \geq c\delta^C.\]
\end{proposition}

It might be conjectured that the first of these bounds can be improved to $\dim(\S) \leq C\log(1/\delta)$ and the second to $|\S| \geq c\delta^C|A|$. This might be called a \emph{Weak Polynomial Freiman-Ruzsa Conjecture} by analogy with \cite{green-fin-field}.

The final result of this section is the one we shall actually use in the sequel. It has the same form as Proposition \ref{weak-bsg}, but in place of the condition that there are many additive quadruples we impose a condition which may appear rather strange at first sight, but is designed specifically with the application we have in mind in the next section.

If $A = \{a_1,\dots,a_k\}$ is a subset of an abelian group $G$ then we say that $A$ is \emph{dissociated} if the only solution to $\eps_1 a_1 + \dots + \eps_k a_k = 0$ with $\eps_i \in \{-1,0,1\}$ is the trivial solution in which $\eps_i = 0$ for all $i$. Recall also that $\langle A \rangle$ denotes the set of all sums $\eps_1 a_1 + \dots + \eps_k a_k$ with $\eps_i \in \{-1,0,1\}$ for all $i$.

\begin{definition}[Arithmetic connectedness]
Suppose that $A \subseteq G$ is a set with $0 \notin A$ and that $m \geq 1$ is an integer. We say that $A$ is $m$-arithmetically connected if, for any set $A' \subseteq A$ with $|A'| = m$ we have either
\begin{enumerate}
\item $A'$ is not dissociated or
\item $A'$ is dissociated, and there is some $x \in A \setminus A'$ with $x \in \langle A' \rangle$.
\end{enumerate}
\end{definition}

\begin{proposition}[Arithmetic connectedness and Bourgain systems]\label{to-use-sec5}
Suppose that $m \geq 1$ is an integer, and that a set $A$ in some abelian group $G$ is $m$-arithmetically connected. Suppose that $0 \notin A$. Then there is a regular Bourgain system $\S$ satisfying
\[ \dim(\S) \leq e^{Cm};\]
\[ |\S| \geq e^{-e^{Cm}}|A|\]
and
\[ \Vert \psi_{\S}1_A \Vert_{\infty} \geq e^{-Cm}.\]
\end{proposition}
\proof By Proposition \ref{weak-bsg}, it suffices to prove that an $m$-arithmetically connected set $A$ has at least $e^{-Cm}|A|^3$ additive quadruples. If $|A| < m^2$ this result is trivial, so we stipulate that $|A| \geq m^2$. Pick any $m$-tuple $(a_1,...,a_m)$ of distinct
elements of $A$. With the stipulated lower bound on $|A|$, there
are at least $|A|^m/2$ such $m$-tuples. We know that either the
vectors $a_1,...,a_m$ are not dissociated, or else there is a
further $a' \in A$ such that $a'$ lies in the linear span of the
$a_i$. In either situation there is some non-trivial linear
relation
\begin{equation*}
\lambda_1 a_1 + \dots + \lambda_m a_m + \lambda' a' = 0
\end{equation*}
where $\vec{\lambda} := (\lambda_1,\dots,\lambda_m,\lambda')$ has
elements in $\{-1,0,1\}$ and, since $0 \notin A$ and the $a_i$s
(and $a'$) are distinct, at least three of the components of
$\vec{\lambda}$ are nonzero. By the pigeonhole principle, it
follows that there is some $\vec{\lambda}$ such that the linear
equation
\begin{equation*}
 \lambda_1 x_1 + \dots+ \lambda_m x_m + \lambda' x' = 0
\end{equation*}
has at least $\frac{1}{2 \cdot 3^{m+1}}|A|^m$ solutions with $x_1,...,x_m,x' \in
A$. Removing the zero coefficients, we may thus assert that there
are some non-negative integers $r_1,r_2$, $3 \leq r_1+r_2
\leq m + 1$, such that the equation
\begin{equation*}
x_1 + \dots + x_{r_1} - y_1- \dots - y_{r_2}=0
\end{equation*}
has at least $\frac{1}{6m^2 3^m}|A|^{r_1+r_2-1} \geq e^{-Cm}|A|^{r_1 + r_2 - 1}$ solutions with
$x_1,\dots,x_{r_1},y_1,\dots,y_{r_2}\in A$. Note that this is a
strong structural statement about $A$, since the maximum possible
number of solutions to such an equation is $|A|^{r_1+r_2-1}$.

We may deduce directly from this the claim that there are at least $e^{-C'm}|A|^3$ additive quadruples in $A$. To do this observe that what we have shown may be recast in the form
\[ 1_A \ast \dots \ast 1_A \ast 1_{-A} \ast \dots \ast 1_{-A}(0) \geq e^{-Cm}\Vert 1_A \Vert_1^{r_1 + r_2 - 1},\]
where there are $r_1$ copies of $1_A$ and $r_2$ copies of $1_{-A}$.

Writing this in terms of the Fourier transform gives
\[ \Vert \widehat{1}_A \Vert_{r_1 + r_2}^{r_1 + r_2} \geq \sum_{\gamma} \widehat{1}_A(\gamma)^{r_1}\widehat{1}_{A}(\overline{\gamma})^{r_2} \geq e^{-Cm}\Vert 1_A \Vert_1^{r_1 + r_2 - 1}.\]
By H\"older's inequality this implies that
\begin{equation}\label{use-soon-10} \Vert \widehat{1}_A \Vert_4^2 \Vert \widehat{1}_A \Vert_{2r_1 + 2r_2 - 4}^{r_1 + r_2 - 2} \geq e^{-Cm} \Vert 1_A \Vert_1^{r_1 + r_2 - 1}.\end{equation}
However if $k$ is an integer then $\Vert \widehat{1}_A \Vert_{2k}^{2k}$ is $|G|^{1-2k}$ times the number of solutions to $a_1 + \dots + a_k = a'_1 + \dots + a'_k$ with $a_i, a'_i \in A$, and this latter quantity is clearly at most $|A|^{2k-1}$. Thus 
\[ \Vert \widehat{1}_A \Vert_{2k} \leq \Vert 1_A \Vert^{1 - 1/2k}_1.\] (In fact, the same is true if $2k \geq 2$ is any real number, by the Hausdorff-Young inequality and the fact that $f$ is bounded by $1$.) Setting $k = r_1 + r_2 - 2$ and substituting into \eqref{use-soon-10}, we immediately obtain
\[ \Vert \widehat{1}_A \Vert_4^4 \geq e^{-2Cm} \Vert 1_A \Vert_1^3,\] which is equivalent to the result we claimed about the number of additive quadruples in $A$.\endproof

\section{Concentration on a Bourgain system}

\begin{proposition}\label{prop7.1}
Suppose that $f : G \rightarrow \R$ has $\Vert f \Vert_A \leq M$, $M \geq 1/2$, and $d(f,\Z) \leq e^{-CM^4}$. Then there is a regular Bourgain system $\S$ with 
\[ \dim(\S) \leq e^{CM^4},\]
\[ \mu(\S) \geq e^{-e^{CM^4}}\Vert f_{\Z} \Vert_1\] and
\[ \Vert \psi_{\S}f \Vert_{\infty} \geq e^{-CM^4}.\]
\end{proposition}
\proof We first obtain a similar result with $g := f^2$ replacing $f$. This function, of course, has the advantage of being non-negative. Note that $\Vert g \Vert_A \leq M^2$, and also that
\[ \Vert g - f_{\Z}^2 \Vert_{\infty} \leq \Vert f - f_{\Z}\Vert_{\infty} \Vert f + f_{\Z} \Vert_{\infty} \leq d(f,\Z) (2 \Vert f \Vert_A + 1) \leq 4Md(f,\Z),\]
and so $d(g,\Z) \leq 4Md(f,\Z)$.

Write $A := \Supp(g_{\Z})$, and $m := \lceil 50M^4\rceil$. If $A= G$ the result is trivial; otherwise, by subjecting $f$ to a suitable translation we may assume without loss of generality that $0 \notin A$. We claim that $A$ is $m$-arithmetically connected. If this is not the case then there are dissociated elements $a_1,\dots,a_m \in A$ such that there is no further $x \in A$ lying in the span $\langle a_1,\dots,a_m\rangle$. Consider the function $p(x)$ defined using its Fourier transform by the Riesz product
\[ \widehat{p}(\gamma) := \prod_{i = 1}^m (1 + \frac{1}{2}(\gamma(a_i) + \overline{\gamma}(a_i))).\]
It is easy to check that $p$ enjoys the standard properties of Riesz products, namely that $\widehat{p}$ is real and nonnegative and that $\Vert p \Vert_A = \sum_{\gamma} \widehat{p}(\gamma) = 1$, and that $\Supp(p) \subseteq \langle a_1,\dots,a_m\rangle$.

Thus we have
\[ \Vert g p \Vert_A \leq \Vert g \Vert_A \Vert p \Vert_A \leq M^2\]
and 
\[ \Vert (g - g_{\Z})p \Vert_A \leq \sum_{x \in \langle a_1,\dots,a_m\rangle} \Vert (g - g_{\Z}) p1_x\Vert_A \leq 3^m \Vert g - g_{\Z}\Vert_{\infty} \leq 4 M3^m d(f,\Z).\]
Thus, since $d(f,\Z)$ is so small, we have
\[ \Vert g_{\Z} p \Vert_A \leq 2M^2.\]
Now since $A \cap \langle a_1 ,\dots, a_m \rangle = \{a_1,\dots,a_m\}$ we have
\[ g_{\Z}p(x) = \sum_{i = 1}^m g_{\Z}(a_i)p(a_i) 1_{a_i}(x).\]
Noting that
\[ p(a_i) = \sum_{\vec{\eps}: \eps_1a_1 + \dots + \eps_m a_m = a_i} 2^{-\sum_j|\eps_j|} \geq \frac{1}{2},\]
we see that
\[ \Vert \widehat{g_{\Z} p} \Vert_2^2 = \Vert g_{\Z} p \Vert_2^2 \geq \frac{1}{4|G|}\sum_{i = 1}^m |g_{\Z}(a_i)|^2 \geq \frac{m}{4|G|}\] and
\begin{align*} \Vert \widehat{g_{\Z} p} \Vert_4^4 &= \frac{1}{|G|^3} \sum_{\substack{i_1,i_2,i_3,i_4 \\ a_{i_1} + a_{i_2} = a_{i_3} + a_{i_4}}} |g_{\Z}(a_i) p(a_i)|^4 \\ & \leq \frac{3}{|G|^3} \big( \sum_{i = 1}^m |g_{\Z}(a_i) p(a_i)|^2  \big)^2 \leq \frac{3}{|G|} \Vert \widehat{g_{\Z} p} \Vert_2^4,
\end{align*}
the middle inequality following from the fact that $a_{i_1} + a_{i_2} = a_{i_3} + a_{i_4}$ only if $i_1 = i_3, i_2 = i_4$ or $i_1 = i_4, i_2 = i_3$ or $i_1 = i_2, i_3 = i_4$. From H\"older's inequality we thus obtain
\[ \Vert g_{\Z}p \Vert_A \geq \frac{\Vert \widehat{g_{\Z}p} \Vert_2^3}{\Vert \widehat{g_{\Z}p} \Vert_4^2} \geq \sqrt{ \frac{|G|}{3}} \Vert \widehat{g_{\Z}p} \Vert_2 \geq \sqrt{\frac{m}{12}}.\]
Recalling our choice of $m$, we see that this contradicts the upper bound $\Vert g_{\Z}p \Vert_A \leq 2M^2$ we obtained earlier.

This proves our claim that $A = \Supp(g_{\Z})$ is $50M^4$-arithmetically connected. It follows from Proposition \ref{to-use-sec5} that there is a regular Bourgain system $\S = (X_{\rho})_{\rho \in [0,4]}$ with 
\begin{equation}\label{eq6a} \dim(\S) \leq e^{CM^4},\end{equation}
\[ |\S| \geq e^{-e^{CM^4}}|A|\] and
\[ \Vert \psi_{\S} 1_A \Vert_{\infty} \geq e^{-CM^4}.\]
Since $\Vert f \Vert_{\infty} \leq M$, the second of these implies that
\begin{equation}\label{eq6b} \mu(\S) \geq e^{-e^{C'M^4}}\Vert f_{\Z} \Vert_1.\end{equation}
Since $g(x) \geq 1_A(x)/2$, the last of these implies that
\begin{equation}\label{g-corr} \Vert \psi_{\S} g \Vert_{\infty} \geq e^{-CM^4}.\end{equation}

It remains to convert these facts about $g$ to the required facts about $f$. The corresponding argument in \cite{green-sanders} (which comes near the end of Proposition 5.1) is rather short, but in the setting of a general Bourgain system we must work a little harder.

We have proved \eqref{g-corr}, which implies the existence of an $x$ such that 
\[ |\E_y f(y)^2\beta_1(x-y)| = \delta,\] for some $\delta \geq e^{-CM^4}$.
Writing this in terms of the Fourier transform we have
\[ |\sum_{\gamma} (f(\cdot) \beta_1(x - \cdot))^{\wedge}(\gamma) \widehat{f}(\gamma)| = \delta\]
which, since $\Vert f \Vert_A \leq M$, implies that there is a $\gamma \in \widehat{G}$ such that
\[ |(f (\cdot) \beta_1(x - \cdot))^{\wedge}(\gamma)| \geq \delta/M,\]
or in other words
\begin{equation}\label{eq7.7} |\E_y f(y) \beta_1(y-x) \gamma(y)| \geq \delta/M.\end{equation}
Now define \[\S' := \lambda \big(\frac{\delta}{80dM^2}\S \cap \Bohr_{\delta/8M^2}(\{\gamma\})\big),\] where as usual $\lambda \in [1/2,1]$ is chosen so that $\S'$ is regular. Now since $\Supp(\beta'_1) \subseteq X_{\delta/40dM^2}$ we have from Lemma \ref{lemma3.2} that
\[ \Vert \beta_1 \ast \beta'_1 - \beta_1 \Vert_1 = \E_x \big|\E_y \beta'_1(y) (\beta_1(x-y) - \beta_1(x))  \big| \leq \delta/2M^2.\]
Since $\Vert f \Vert_{\infty} \leq \Vert f \Vert_A \leq M$ we may introduce an averaging over $\beta'_1$ into \eqref{eq7.7}, obtaining
\begin{equation}\label{eq7.8} |\E_y f(y)\gamma(y) \E_t \beta_1(y + t - x)\beta'_1(t) | \geq \delta/2M.\end{equation}
Now if $t \in \Supp(\beta'_1)$ then by construction we have
\[ |1 - \gamma(t)| \leq \delta/4M^2,\] and so from \eqref{eq7.8} and the fact that $\Vert f \Vert_{\infty} \leq M$ we see that
\[ |\E_y f(y)\gamma(y + t) \E_t \beta_1(y + t - x) \beta'_1(t)| \geq \delta/4M.\]
Changing variables by writing $z := y + t - x$ and noting that $|\gamma(x)| = 1$, this may be written
\[ |\E_z \beta_1(z)\gamma(z) \E_y f(y)\beta'_1(z + x - y)| \geq \delta/4M,\]
which immediately implies that
\[ \Vert \psi_{\S'}f \Vert_{\infty} \geq \delta/4M \geq e^{-C'M^4}.\]
To conclude the argument we must show that $\S'$ is subject to the same bounds \eqref{eq6a}, \eqref{eq6b} (possibly with different constants $C$). This follows easily from Lemma \ref{join-lem} and the bounds of Lemma \ref{bohr-facts}.\endproof

\section{The inductive step}

Our remaining task is to prove Lemma \ref{ind-step}, the inductive step which drives the proof of Theorem \ref{mainthm-finite}. We recall the statement of that lemma now for the reader's convenience.

\begin{inductive}[Inductive Step]
Suppose that $f : G \rightarrow \R$ has $\Vert f \Vert_A \leq M$, where we take $M \geq 1$, and that $d(f,\Z) \leq e^{-C_1M^4}$. Set $\eps := e^{-C_0M^4}$, for some constant $C_0$.
Then we may write $f = f_1 + f_2$, where 
\begin{enumerate}
\item either $\Vert f_1 \Vert_A \leq \Vert f \Vert_A - 1/2$ or else $(f_1)_{\Z}$ may be written as $\sum_{j = 1}^L \pm 1_{x_j + H}$, where $H$ is a subgroup of $G$ and $L \leq e^{e^{C'(C_0)M^4}}$;
\item $\Vert f_2 \Vert_A \leq \Vert f \Vert_A - \frac{1}{2}$ and 
\item $d(f_1,\Z) \leq d(f,\Z) + \eps$ and $d(f_2,\Z) \leq 2d(f,\Z) + \eps$.
\end{enumerate}
\end{inductive}
\begin{remark} Note carefully the factor $2$ in the bound for $d(f_2,\Z)$; this is one important reason for the weakness of our bounds in Theorem \ref{mainthm-finite}.\end{remark}

\proof 
Applying Proposition \ref{prop7.1} to $f$ we obtain a regular Bourgain system $\S$ with
\[ d = \dim(\S) \leq e^{CM^4},\]
\[ \mu( \S ) \geq e^{-e^{CM^4}} \Vert f_{\Z} \Vert_1\]
and 
\[ \Vert \psi_{\S} f \Vert_{\infty} \geq 3e^{-CM^4}.\]
We were given that $\eps = e^{-C_0M^4}$ in the statement of the proposition. Without loss of generality (by increasing $C_0$) we may assume that $C_0$ is much larger than the constant $C$ in the bounds just given.

Applying Proposition \ref{prop4.2} with this value of $\eps$ we obtain a regular Bourgain system $\S' = (X'_{\rho})_{\rho \in [0,4]}$ such that
\[ \dim(\S') \leq e^{C'M^4},\]
\begin{equation}\label{eq8.5} \mu( \S' ) \geq e^{-e^{C'M^4}} \Vert f_{\Z} \Vert_1\end{equation}
and 
\begin{equation}\label{eq8.6} \Vert \psi_{\S'} f \Vert_{\infty} \geq \Vert \psi_{\S} f \Vert_{\infty} - \eps > 2e^{-CM^4},\end{equation}
and with the additional property that 
\begin{equation}\label{psi-s-int} d(\psi_{\S'}f,\Z) \leq d(f,\Z) + \eps.\end{equation}
Here, $C' = C'(C_0)$ depends only on $C_0$.
We define $f_1 := \psi_{\S'}f$ and $f_2 := f - \psi_{\S'}f$. Thus we immediately see that $d(f_1,\Z) \leq d(f,\Z) + \eps$, which is one part of (iii), and the other inequality
$d(f_2,\Z) \leq 2d(f,\Z) + \eps$ follows immediately from this.

Note also that \eqref{eq8.6} and the fact that $d(f_1,\Z) \leq d(f,\Z) + \eps < 2e^{-CM^4}$ implies that $(f_1)_{\Z}$ is not identically zero, and therefore $\Vert f_1 \Vert_{\infty} \geq \Vert (f_1)_{\Z} \Vert_{\infty} - \eps > 1/2$, and hence $\Vert f_1 \Vert_A \geq \frac{1}{2}$. It follows that $\Vert f_2 \Vert_A \leq \Vert f \Vert_A - \frac{1}{2}$, as we were required to prove.

It remains to deal with the possibility that $\Vert f_1 \Vert_A > \Vert f \Vert_A - \frac{1}{2}$. If this is so then $\Vert f_2 \Vert_A < 1/2$, and thus $\Vert f_2 \Vert_{\infty} < 1/2$. It follows that $(f_2)_{\Z} = 0$, and hence $f_{\Z} = (\psi_{\S'}f)_{\Z}$.

Now suppose that $x - x' \in X'_{\eps/20dM}$. Then from Lemma \ref{lemma3.3} and \eqref{psi-s-int} we see that
\begin{align*} |f_{\Z}(x) - f_{\Z}(x')| &= |(\psi_{\S'}f)_{\Z}(x) - (\psi_{\S'}f)_{\Z}(x')| \\ & \leq |\psi_{\S'}f(x) - \psi_{\S'}f(x')| + 2(d(f,\Z) + \eps) \leq 2d(f,\Z) + 3\eps < 1/10.\end{align*}
We are forced to conclude that $f_{\Z}(x) = f_{\Z}(x')$. It follows immediately that $f_{\Z}$ is constant on cosets of the subgroup $H := \langle X_{\eps/20dM}\rangle$, and so
\[ f_{\Z} = \sum_{j = 1}^L \pm 1_{x_j + H}\] for some $j = 1,\dots,L$, where we may take
\[ L \leq M\Vert f_{\Z} \Vert_1/\Vert 1_H \Vert_1.\]
 Let us note from property \textsc{bs5} of Bourgain systems that 
\[ \Vert 1_H \Vert_1  \geq \frac{|X_{\eps/20dM}|}{|G|} \geq \big( \frac{\eps}{40dM} \big)^d \mu(\S').\]
It follows from this and \eqref{eq8.5} that 
\[ L \leq e^{e^{C''M^4}},\] where $C'' = C''(C_0)$ depends only on $C_0$, as required.\endproof

\section{Possible improvements}

As it stands, our argument ``loses an exponential'' in two places. First of all the ``almost integer'' parameter $d(f,\Z)$ must not be allowed to blow up during the iteration leading to the proof of Theorem \ref{mainthm-finite}. This requires it to be exponentially small in $M$ at the beginning of the argument. This parameter then gets exponentiated again in any application of Proposition \ref{prop4.2}.

We note that our proof in fact yields a version of Theorem \ref{mainthm-finite} for functions $f$ satisfying $d(f,\Z) \leq e^{-CM^4}$, rather than simply $d(f,\Z) = 0$. Such functions are within $d(f,\Z)$ of a sum $\sum_{j = 1}^L \pm 1_{x_j + H_j}$. Now an example of M\'ela \cite{mela} shows that, for such a theorem to hold, $d(f,\Z)$ must be exponentially small in $M$. It seems then that, as long as our proof technique also establishes this more general result, the bound we can hope to obtain is seriously restricted.

The function $f := 1_{\{1,\dots,N\}}$ on the group $G = \Z$ has $\Vert f \Vert_A \sim \log N$, yet it cannot be written as the $\pm$-sum of fewer than $N$ cosets in $\Z$. This shows that the bound of Theorem \ref{mainthm} cannot be improved beyond $L \leq e^{C\Vert \mu \Vert}$ in general. It may be that this represents the correct bound. Note that this would immediately provide another proof of the Littlewood conjecture, to add to the famous papers of Konyagin\cite{konyagin} and McGehee, Pigno and Smith \cite{mps}. Recall that this conjecture was the following statement:  if $A \subseteq \Z$ is a set of size $N$ then

\[ \int^1_0 \big| \sum_{ a \in A} e^{2\pi i \theta a} \big|  \, d\theta \gg \log N.\]

Our results imply the weaker inequality
\[ \int^1_0 \big| \sum_{ a \in A} e^{2\pi i \theta a} \big|  \, d\theta \gg (\log \log N)^{1/4},\]
easily the weakest bound ever obtained in the direction of the Littlewood conjecture!

\appendix

\section{Reduction of Theorem \ref{mainthm} to the finite case}\label{lcag-finite}

 Throughout the section $G$ will be a locally compact abelian group.

In this section we deduce Theorem \ref{mainthm} from Theorem \ref{mainthm-finite}. As we remarked, this section is independent of the rest of the paper. It is also not self-contained, and in particular we assume the (qualitative) idempotent theorem. The reader may safely think of the case $G = \T^d$, $\widehat{G} = \Z^d$, which captures the essence of the argument and may be thought of in quite concrete terms. 

We begin by proving the following special case of Theorem \ref{mainthm}.

\begin{proposition}[The finite case]\label{finite-gen-prop}
Suppose that $G$ is compact, and that $\mu \in \M(G)$ satisfies $\Vert \mu \Vert \leq M$ and has the form
\[ \widehat{\mu} = \sum_{j = 1}^K \pm 1_{\gamma_j}.\]
Then we may write
\[ \widehat{\mu} = \sum_{l = 1}^L \pm 1_{\gamma'_l + \Gamma_l},\]
where $L \leq e^{e^{CM^4}}$, each $\Gamma_l$ is a subgroup of $\widehat{G}$ and there are at most $\Vert \mu \Vert + \frac{1}{100M}$ different groups $\Gamma_l$.
\end{proposition}
\proof  We may suppose that $\widehat{G} = \langle \gamma_1,\dots,\gamma_K\rangle$. By the structure theorem for finitely-generated abelian groups, $\widehat{G}$ is isomorphic to $\widehat{H} \times \Z^d$, where $H$ is finite. We may now explicitly describe the characters $\gamma_j$ as
\[ \gamma_j = (\omega^{(j)}, r_1^{(j)}, \dots, r_d^{(j)}),\]
where $\omega^{(j)} \in \widehat{H}$ and the $r_i^{(j)}$ are integers. 

Let $N = N(K)$ be an enormous prime and define the measure $\widetilde{\mu}$ on $H \times (\Z/N\Z)^d$ by
\[ \widetilde{\mu}(h,x_1,\dots, x_d) := \sum_{j = 1}^K \omega^{(j)}(h) e\big(\frac{r_1^{(j)} x_1 + \dots + r_d^{(j)}x_d}{N}  \big).\]

It is clear that
\begin{align*}\lim_{N \rightarrow \infty} \Vert \widetilde{\mu} \Vert &= \lim_{N \rightarrow \infty} \E_{h \in H, x_1,\dots, x_d \in \Z/N\Z} \big| \sum_{j = 1}^K \omega^{(j)}(h) e\big( \frac{r_1^{(j)} x_1 + \dots + r_d^{(j)}x_d}{N} \big)  \big|\\ &= \E_h \int_{\theta_1,\dots, \theta_d \in \T^d} \big| \sum_{j = 1}^K \omega^{(j)}(h) e(r^{(j)}_1 \theta_1 + \dots + r^{(j)}_d \theta_d) \big|\, d\theta_1 \dots d\theta_d = \Vert \mu \Vert.\end{align*}
Taking $N$ large enough, we may assume that $\Vert \widetilde{\mu} \Vert \leq \Vert \mu \Vert + \frac{1}{200M}$.

Theorem \ref{mainthm-finite} now applies\footnote{In fact, this is not quite true. Theorem \ref{mainthm-finite} stated a weaker bound of $\Vert \widetilde{\mu} \Vert + \frac{1}{100}$ on the number of distinct $\Gamma_j$, but that was done only for clarity and it is clear that a trivial modification of the proof gives this somewhat stronger bound whilst preserving the upper bound on $L$.} with $f := \widehat{\widetilde{\mu}}$ to show that 
\begin{equation}\label{eq4.5} \widehat{\widetilde{\mu}} = \sum_{l = 1}^L \pm 1_{\gamma'_l + \Gamma_l},\end{equation}
where the $\Gamma_l$ are subgroups of $\widehat{H} \times (\Z/N\Z)^d$, $L \leq e^{e^{CM^4}}$, and we may assume the number of distinct $\Gamma_l$ is at most
\[ \Vert \widetilde{\mu} \Vert + \frac{1}{200M} \leq \Vert \mu \Vert + \frac{1}{100M}.\]

Let $\pi : \widehat{H} \times (\Z/N\Z)^d \rightarrow \widehat{H}$ be the obvious projection, and for each group $\Gamma_l$ appearing in the decomposition \eqref{eq4.5} consider the subgroup $\Gamma_l \cap \ker\pi$. Suppose without loss of generality that $|\Gamma_1 \cap \ker\pi| \geq |\Gamma_l \cap \ker\pi|$ for $l = 1,\dots,L$, and also that $\Gamma_1 \cap \ker\pi = \Gamma_l \cap \ker\pi$ precisely for $l = 1,\dots,m$.

We clearly have $|\Gamma_1 \cap \ker\pi| = N^{d'}$ for some $d'$, $0 \leq d' \leq d$. Suppose that $d' \geq 1$. The portion $\sum_{l = 1}^m \pm 1_{\gamma'_l + \Gamma_l}$ is constant on cosets of $\Gamma_1 \cap \ker\pi$.  If this portion is zero then we may simply delete it from \eqref{eq4.5}. Suppose, then, that it is not zero. We note that if $l \in \{1,\dots,m\}$ and  $k \notin \{1,\dots,m\}$ then 
\[ |(\gamma'_l + \Gamma_l) \cap (\gamma'_k + \Gamma_k)| \leq |H|N^{d' - 1}.\]
We therefore have the estimate
\[ \Vert \sum_{l = 1}^L \pm 1_{\gamma'_l + \Gamma_l} \Vert_1 \geq |\Gamma_1 \cap \ker\pi| - |H|N^{d' - 1}m(L-m) \geq N^{d'} - |H|N^{d'-1}m(L - m).\] If $N$ is chosen large enough this gives us
\[  \Vert \sum_{l = 1}^L \pm 1_{\gamma'_l + \Gamma_l} \Vert_1 > K.\] However by definition $\widehat{\widetilde{\mu}}$ is the characteristic function of a set of $K$ characters, and so we have a contradiction.

It follows that we may assume $d' = 0$, in which case all of the subgroups $\Gamma_l$ appearing in \eqref{eq4.5} are simply subgroups of $\widehat{H}$. The decomposition \eqref{eq4.5} may then be regarded as a decomposition of $\widehat{\mu}$ as well, and we have proved our result.\endproof

It remains to reduce to the finite case covered in this last proposition. Let $\mu \in \M(G)$ be an arbitrary idempotent measure, and write $M := \Vert \mu \Vert$. The idempotent theorem implies that $\widehat{\mu}$ may be written as a finite combination
\[ \widehat{\mu} = \sum_{k \in E} \pm 1_{\gamma_k + \Gamma_k}.\]
We say that two open subgroups $\Gamma_k,\Gamma_l$ are \emph{commensurable} if $|\Gamma_k : \Gamma_k \cap \Gamma_l|, |\Gamma_l : \Gamma_k \cap \Gamma_l| < \infty$; this is an equivalence relation. Split $E$ into a disjoint union $E_1 \cup \dots \cup E_J$ of equivalence classes.
We have 
\[ \widehat{\mu} = \sum_{k = 1}^K \pm 1_{\gamma_k + \Gamma_k} = \sum_{j = 1}^J \sum_{k \in E_j} \pm 1_{\gamma_k + \Gamma_k}.\]
Writing $\Omega_j := \bigcap_{k \in E_j} \Gamma_k$ and noting that $|\Gamma_k : \Omega_j| < \infty$ for $k \in E_j$, we may write this in the form
\[ \widehat{\mu} = \sum_{j = 1}^J \sum_{k \in E'_j} \pm 1_{\gamma_k + \Omega_j},\]
where the index sets $E'_j$ are still finite and the open subgroups $\Omega_j$, $j = 1,\dots,J$, are mutually incommensurable.

It follows that
\[ \mu = \sum_{j = 1}^J \mu_j,\]
where
\[ \mu_j(x) := \sum_{k \in E'_j} \pm \gamma_k(x) \mu_{H_j}(x).\]
Here, $\mu_{H_j}$ is the Haar measure on the compact group $H_j := \Omega_j^{\perp}$. The incommensurability of the $H_j$ implies that
\[ \Vert \mu \Vert = \sum_{j = 1}^J \Vert \mu_j \Vert.\]
The measures $\mu_j$ need not be idempotent, but their Fourier-Stieltjes transforms are, by construction, integer-valued. We may of course suppose that no $\mu_j$ is zero, and hence we have the upper bound $J \leq \Vert \mu \Vert = M$.

We may now simply apply Proposition \ref{finite-gen-prop} to $\mu_j$ which, when regarded as a measure on $H_j$, is of the form covered by that proposition.

Once this is done, we shall have written $\widehat{\mu}$ as a sum 
\[ \widehat{\mu} = \sum_{q = 1}^Q \pm 1_{\gamma_q + \Gamma_q} \]where
\[ Q \leq J e^{e^{CM^4}} \leq e^{e^{C'M^4}}\]
and the number of distinct $\Gamma_q$ is bounded by
\[ \sum_{j = 1}^J \big( \Vert \mu_j \Vert + \frac{1}{100M}\big) \leq \Vert \mu \Vert + \frac{1}{100}.\]
This concludes the proof of Theorem \ref{mainthm}.\endproof


\begin{thebibliography}{99}


\bibitem{bourgain} J.~Bourgain, \emph{On triples in arithmetic progression,} Geom. Fuct. Anal.  \textbf{9}  (1999),  no. 5, 968--984. 

\bibitem{bog} N.~N.~Bogolyubov, \emph{Sur quelques propri\'et\'es arithm\'etiques des presque-p\'eriodes,} Ann. Chaire Math. Phys. Kiev \textbf{4} (1939), 185--194.

\bibitem{chang-freiman} M.-C.~Chang, \emph{A polynomial bound in Freiman's theorem,}  Duke Math. J.  \textbf{113}  (2002),  no. 3, 399--419. 
 

\bibitem{cohen} P.~J.~Cohen, \emph{On a conjecture of Littlewood and Idempotent measures,} Amer. J. Math. \textbf{82}, no. 2, (1960), 191--212.

\bibitem{croft-eureka} H.~T.~Croft, \emph{Some problems,} Eureka, Cambridge University 1968.

\bibitem{gowers-4aps} W.~T.~Gowers, \emph{A new proof of Szemer\'edi's theorem for arithmetic progressions of length four,} Geom. Funct. Anal. \textbf{8} (1998), 529--551.

\bibitem{green-fin-field} B.~J.~Green, \emph{Finite field models in additive combinatorics,} in Surveys in Combinatorics 2005, London Math. Soc. Lecture Notes \textbf{327}, 1--27.

\bibitem{green-regularity} \bysame, \emph{A Szemer\'edi-type regularity lemma in abelian groups, with applications,}  Geom. Funct. Anal.  \textbf{15}  (2005),  no. 2, 340--376.

\bibitem{green-nonabelian} \bysame, \emph{Notes on sets with small doubling in non-abelian groups,} in preparation.

\bibitem{green-kony} B.~J.~Green and S.~V.~Konyagin, \emph{The Littlewood problem modulo a prime,} to appear in Canadian J. Math.

\bibitem{green-ruzsa} B.~J.~Green and I.~Z.~Ruzsa, \emph{An analogue of Freiman's theorem in an arbitrary abelian group,} to appear in J. London Math. Soc.

\bibitem{green-sanders} B.~J.~Green and T.~Sanders, \emph{Boolean functions with small spectral norm,} preprint available at\\
\texttt{http://www.arxiv.org/abs/math.CA/0605524.}

\bibitem{green-tao-inverseu3} B.~J.~Green and T.~C.~Tao, \emph{An inverse theorem for the Gowers $U^3$-norm, with applications,} to appear in Proc. Edinburgh Math. Soc.


\bibitem{heinonen} J.~Heinonen, \emph{Lectures on analysis on metric spaces,} Springer Universitext 2001.

\bibitem{helson} H.~Helson, \emph{Note on harmonic functions,} Proc. Amer. Math. Soc. \textbf{4} (1953), 686--691.

\bibitem{konyagin} S.~V.~Konyagin, \emph{On the Littlewood problem,} (Russian)
Izv. Akad. Nauk SSSR Ser. Mat. \textbf{45} (1981), no. 2, 243--265, 463. English translation in Math. USSR Izvestija \textbf{18} (1982), no. 2, 205--225.


\bibitem{mps} O.~C.~McGehee, L.~Pigno, and B.~Smith, \emph{Hardy's inequality and the $L\sp{1}$ norm of exponential sums,}  Ann. of Math. (2)  \textbf{113}  (1981), no. 3, 613--618.

\bibitem{mela} J.~-F.~M\'ela, \emph{Mesures $\eps$-idempotentes de norme born\'ee,} Studia Math. \textbf{72} (1982), no. 2, 131--149.

\bibitem{rudin-book} W.~Rudin, \emph{Fourier analysis on groups,} 2nd Ed., Wiley 1990.

\bibitem{rudin-paper} \bysame, \emph{Idempotent measures on abelian groups,} Pacific J. Math. \textbf{9} (1959) 195--209.


\bibitem{sanders1} T.~Sanders, \emph{An application of a local version of Chang's theorem,} preprint. Available at\\
\texttt{http://www.arxiv.org/abs/math.CA/0607668}.

\bibitem{sanders2} \bysame, \emph{The Littlewood-Gowers problem,} preprint. Available at\\
\texttt{http://www.arxiv.org/abs/math.CA/0605523}.

\bibitem{shkredov} I.~D.~Shkredov, \emph{On a generalization of Szemer\'edi's theorem,} preprint. Available at\\
\texttt{http://www.arxiv.org/abs/math.NT/0503639}.

\bibitem{tao:bourgain} T.~C.~Tao, \emph{The Roth-Bourgain theorem,} ``Short story'' available at\\
 \texttt{http://www.math.ucla.edu/$\widetilde{\;}$tao}.

\bibitem{tao-vu-book} T.~C.~Tao and V.~H.~Vu, \emph{Additive combinatorics,} Cambridge Studies in Advanced Mathematics \textbf{105}, CUP 2006.

\end{thebibliography}
\end{document}